\documentstyle[11pt,a4wide]{article}
\input amssym.def
\input amssym

\font \tenmsb=msbm10 scaled \magstep 1
\font \sevenmsb=msbm7 scaled \magstep 1
\font \fivemsb=msbm5 scaled \magstep 1
\textfont \msbfam=\tenmsb
\scriptfont \msbfam=\sevenmsb
\scriptscriptfont \msbfam=\fivemsb

\font\teneufm=eufm10 scaled \magstep 1
\font\seveneufm=eufm7 scaled \magstep 1
\font\fiveeufm=eufm5 scaled \magstep 1
\textfont\eufmfam=\teneufm
\scriptfont\eufmfam=\seveneufm
\scriptscriptfont\eufmfam=\fiveeufm

\title{\bf QUANTUM MATRIX BALL: THE WEIGHTED BERGMAN KERNELS}

\author{\sl D. Shklyarov \and \sl S. Sinel'shchikov \and \sl L. Vaksman}

\date{\tt Institute for Low Temperature Physics \& Engineering\\
National Academy of Sciences of Ukraine}

\begin{document}

\maketitle

\bigskip

 The Cartan domains are among the important subjects in many problems of
representation theory and mathematical physics \cite{BE, Hur}. The methods
of quantum groups theory \cite{CP} were used in \cite{SV} to produce
q-analogues of Cartan domains, in particular, q-analogues of balls in the
spaces of complex matrices.

 The point of this work is to consider those quantum matrix balls and the
associated Hilbert spaces of 'functions'. As a main result, we present an
explicit formula for the weighted Bergman kernel.

 It is implicit everywhere in the sequel that $q \in(0,1)$, $m,n \in{\Bbb
N}$, and $m \le n$.

 We need q-analogues for $*$-algebras ${\rm Pol}({\rm Mat}_{mn})$ of
polynomials on the space ${\rm Mat}_{mn}$ of complex matrices and a
$*$-algebra $D({\Bbb U})$ of smooth finite functions in the matrix ball
${\Bbb U}=\{A \in{\rm Mat}_{mn}|\;\|A \|<1 \}$. We start with forming a
q-analogue for the algebra ${\rm Fun}({\Bbb U})={\rm Pol}({\rm
Mat}_{mn})+D({\Bbb U})$.

 The $*$-algebra ${\rm Fun}({\Bbb U})_q$ is given by its generators $f_0$,
$z_a^\alpha$, $a=1,2,\ldots,n$, $\alpha=1,2,\ldots,m$, and the relations
\begin{equation}\label{z's}\left \{\begin{array}{lcl}
z_a^\alpha z_b^\beta=qz_b^\beta z_a^\alpha &,& a=b \quad \&\quad
\alpha<\beta \qquad{\rm or}\qquad a<b \quad \&\quad \alpha=\beta \\
z_a^\alpha z_b^\beta=z_b^\beta z_a^\alpha &,& \alpha<\beta \quad\&\quad a>b
\\
z_a^\alpha z_b^\beta-z_b^\beta z_a^\alpha=(q-q^{-1})z_a^\beta z_b^\alpha &,&
\alpha<\beta \quad \&\quad a<b
\end{array}\right.
\end{equation}
\begin{equation}\label{zsz}
\left(z_b^\beta \right)^*\cdot z_a^\alpha=q^2 \sum_{a',b'=1}^n
\sum_{\alpha',\beta'=1}^mR(b,a,b',a')R(\beta,\alpha,\beta',\alpha')\cdot
z_{a'}^{\alpha'}\cdot \left(z_{b'}^{\beta'}
\right)^*+(1-q^2)\delta_{ab}\delta^{\alpha \beta},
\end{equation}
\begin{equation}\label{f0zc}\left(z_a^\alpha \right)^*f_0=f_0z_a^\alpha=0,
\end{equation}
$$f_0=f_0^*=f_0^2.$$

 Here $a,b=1,2,\ldots,n$, $\alpha,\beta=1,2,\ldots,m$,
$$R(i,j,i',j')=\left \{
\begin{array}{ccl}q^{-1} &,&i \ne j \quad \& \quad i=i'\quad \& \quad j=j'
\\
1 &,& i=j=i'=j' \\
-(q^{-2}-1) &,& i=j \quad \& \quad i'=j' \quad \& \quad j'>j \\
0 &,& otherwise
\end{array}\right..$$

 In this setting, the $*$-subalgebra ${\rm Pol}({\rm Mat}_{mn})_q
\subset{\rm Fun}({\Bbb U})_q$ generated by $z_a^\alpha$, $a=1,2,\ldots,n$,
$\alpha=1,2,\ldots,m$ is a q-analogue of the $*$-algebra ${\rm Pol}({\rm
Mat}_{mn})$, and the bilateral ideal $D({\Bbb U})_q={\rm Pol}({\rm
Mat}_{mn})_qf_0{\rm Pol}({\rm Mat}_{mn})_q$ is a q-analogue of the
$*$-algebra $D({\Bbb U})$. (The element $f_0$ works here as a
delta-function, as one can see from (\ref{f0zc})).

 To motivate our subsequent constructions, observe that (see \cite{SSV}) the
$*$-algebra $D({\Bbb U})_q$ is a $U_q \frak{su}_{nm}$-module algebra
\cite{CP}. Remind the explicit formula for invariant integral from
\cite{SSV}.

 Consider the representation $T$ of ${\rm Fun}({\Bbb U})_q$ in the space
${\cal H}={\rm Fun}({\Bbb U})_qf_0={\rm Pol}({\rm Mat}_{mn})_qf_0$:
$$T(f) \psi=f \psi,\qquad f \in{\rm Fun}({\Bbb U})_q,\quad \psi \in{\cal
H}.$$ There exists a unique positive scalar product in ${\cal H}$ such that
$(f_0,f_0)=1$, and
$$(T(f) \psi_1,\psi_2)=(\psi_1,T({f^*})\psi_2),\qquad f \in{\rm Fun}({\Bbb
U})_q,\quad \psi_1,\psi_2 \in{\cal H}.$$
One can prove that the $*$-algebra ${\rm Pol}({\rm Mat}_{mn})_q$ admits a
unique up to unitary equivalence {\sl faithful} irreducible
$*$-representation by {\sl bounded} operators in a Hilbert space. This
$*$-representation can be produced via extending the operators $T(f)$, $f
\in{\rm Pol}({\rm Mat}_{mn})_q$, onto the completion of the pre-Hilbert
space ${\cal H}$.

 The invariant integral is of the form (see \cite{SSV}):
\begin{equation}\label{ii}
\int \limits_{{\Bbb U}_q}fd \nu={\rm tr}(T(f)q^{-2
\Gamma(\check{\rho})})),\qquad f \in D({\Bbb U})_q,
\end{equation}
with $\Gamma:{\frak h}\to{\rm End}({\cal H})$ being a subrepresentation
of the natural representation of the Cartan subalgebra ${\frak h} \subset
\frak{sl}_N$ in ${\rm Fun}({\Bbb U})_q$, and $\check{\rho}\in{\frak h}$ the
element of this Cartan subalgebra determined by the half sum of positive
roots $\rho$ under the standard pairing of ${\frak h}$ and ${\frak h}^*$.
(To see that this integral is well defined, observe that the operators
$T(f)$, $f \in D({\Bbb U})_q$, are finite dimensional, and ${\cal H}$ is
decomposable into a sum of weight subspaces associated to non-negative
weights.)

 Our immediate intention is to produce q-analogues of weighted Bergman
spaces. In the case $q=1$ one has
\begin{equation}\label{de}
{\rm det}(1-{\bf zz^*})=1+\sum_{k=1}^m(-1)^k{\bf z}^{\wedge k}{\bf
z}^{*\wedge k},
\end{equation}
with ${\bf z}^{\wedge k}$, ${\bf z}^{*\wedge k}$ being the "exterior powers"
of the matrices ${\bf z}$, ${\bf z}^*$, that is, matrices formed by the
minors of order $k$. The operators $(1-q^2)^{-1/2}T(z_a^\alpha)$,
$(1-q^2)^{-1/2}T((z_a^\alpha)^*)$, $a=1,2,\ldots,n$, $\alpha=1,2,\ldots,m$,
are respectively the q-analogues of creation and annihilation operators. The
"creation operators" are placed in the right hand side of (\ref{de}) to the
left of the "annihilation operators".  This allows one to produce a
q-analogue of the polynomial ${\rm det}(1-{\bf zz^*})$ in a standard way as
follows.

 Let $1 \le \alpha_1<\alpha_2<\ldots<\alpha_k \le m$, $1 \le
a_1<a_2<\ldots<a_k \le n$. Introduce q-analogues of minors for the matrix
${\bf z}$:
$${{\bf z}^{\wedge k \,}}_{\{a_1,a_2,\ldots,a_k \}}
^{\{\alpha_1,\alpha_2,\ldots,\alpha_k \}}=\sum_{s \in
S_k}(-q)^{l(s)}z_{a_1}^{\alpha_{s(1)}}z_{a_2}^{\alpha_{s(2)}}\ldots
z_{a_k}^{\alpha_{s(k)}},$$
with $l(s)={\rm card}\{(i,j)|\;i<j \quad \& \quad s(i)>s(j)\}$ being the
length of the permutation $s$.

 The q-analogue $y \in{\rm Pol}({\rm Mat}_{mn})_q$ for the polynomial ${\rm
det}(1-{\bf zz^*})$ is defined by
$$y=1+\sum_{k=1}^m(-1)^k \sum_{\{J'|\;{\rm card}(J')=k \}}\sum_{\{J''|\;{\rm
card}(J'')=k \}}{{\bf z}^{\wedge k \,}}_{J''}^{J'}\cdot \left({{\bf
z}^{\wedge k \,}}_{J''}^{J'}\right)^*.$$

 Let $\lambda > m+n-1$. Now (\ref{ii}) allows one to define the integral
with weight $y^\lambda$ as follows:
$$\int \limits_{{\Bbb U}_q}fd \nu_ \lambda \stackrel{\rm
def}{=}C(\lambda){\rm tr}(T(f)T(y)^\lambda q^{-2
\Gamma(\check{\rho})}),\qquad f \in D({\Bbb U})_q,$$
with $C(\lambda)=\displaystyle
\prod_{j=0}^{n-1}\prod_{k=0}^{m-1}(1-q^{2(\lambda+1-N)}q^{2(j+k)})$ being a
normalizing multiple that provides $\displaystyle \int \limits_{{\Bbb
U}_q}1d \nu_ \lambda=1.$

 The Hilbert space $L^2(d \nu_ \lambda)_q$ is defined as a completion of the
space $D({\Bbb U})_q$ of finite functions with respect to the norm $\|f \|_
\lambda=\left(\displaystyle \int \limits_{{\Bbb U}_q}f^*fd \nu_ \lambda
\right)^{1/2}$. The closure $L_a^2(d \nu_ \lambda)_q$ in $L^2(d \nu_
\lambda)_q$ of the unital subalgebra ${\Bbb C}[{\rm Mat}_{mn}]_q \subset{\rm
Pol}({\rm Mat}_{mn})_q$ generated by $z_a^\alpha$, $a=1,2,\ldots,n$,
$\alpha=1,2,\ldots,m$, will be called a weighted Bergman space.

 Note that the relations (\ref{z's}) and the algebra ${\Bbb C}[{\rm
Mat}_{mn}]_q$ were considered in many works on quantum groups (see
\cite{CP}). The $*$-algebra ${\rm Pol}({\rm Mat}_{mn})_q$ determined by the
relations (\ref{z's}) and (\ref{zsz}), is a q-analogue of the Weyl algebra.
This becomes plausible after one performs the 'change of variables' as
follows:
$$z_a^\alpha \mapsto(1-q^2)^{-{1/2}}z_a^\alpha;\qquad
(z_a^\alpha)^*\mapsto(1-q^2)^{-{1/2}}(z_a^\alpha)^*.$$

 Consider the orthogonal projection $P_ \lambda$ in $L^2(d \nu_ \lambda)_q$
onto the weighted Bergman space $L_a^2(d \nu_ \lambda)_q$. It is possible to
show that $P_ \lambda$ could be written as an integral operator (see
\cite{V})
\begin{equation}\label{pl} P_ \lambda f=\int \limits_{{\Bbb
U}_q}K_ \lambda({\bf z},\hbox{\boldmath $\zeta$}^*)f(\hbox{\boldmath
$\zeta$})d \nu_ \lambda(\hbox{\boldmath $\zeta$}),\qquad f \in D({\Bbb
U})_q.  \end{equation}

 Our intention is to introduce the algebra ${\Bbb C}[[{\rm Mat}_{mn}\times
\overline{\rm Mat}_{mn}]]_q$ of kernels of integral operators and to
determine an explicit form of the weighted Bergman kernel $K_ \lambda
\in{\Bbb C}[[{\rm Mat}_{mn}\times \overline{\rm Mat}_{mn}]]_q$ involved in
(\ref{pl}).

 Introduce the notation
\begin{equation}\label{pk}
\Bbbk_i=\sum_{J'\subset \{1,2,\ldots,m \}\atop{\rm
card}(J')=i}\sum_{J''\subset \{1,2,\ldots,n \}\atop{\rm
card}(J'')=i}{z^{\wedge i \,}}_{J''}^{J'}\otimes \left({z^{\wedge i
\,}}_{J''}^{J'}\right)^*.
\end{equation}

 Let ${\Bbb C}[\overline{\rm Mat}_{mn}]_q \subset{\rm Pol}({\rm Mat}_{mn})_q$
be the unital subalgebra generated by $(z_a^\alpha)^*$, $a=1,2,\ldots,n$,
$\alpha=1,2,\ldots,m$, and ${\Bbb C}[{\rm Mat}_{mn}]_q^{\rm op}$ the algebra
which differs from ${\Bbb C}[{\rm Mat}_{mn}]_q$ by a replacement of its
multiplication law to the opposite one (this replacement is motivated in
\cite{V}). The tensor product algebra ${\Bbb C}[{\rm Mat}_{mn}]_q^{\rm
op}\otimes{\Bbb C}[\overline{\rm Mat}_{mn}]_q$ will be called an algebra of
polynomial kernels. It is possible to show that in this algebra $\Bbbk_i
\Bbbk_j=\Bbbk_j \Bbbk_i$ for all $i,j=1,2,\ldots,m$.

 We follow \cite{SV,SSV} in equipping ${\rm Pol}({\rm Mat}_{mn})_q$ with a
${\Bbb Z}$-grading: ${\rm deg}(z_a^\alpha)=1$, ${\rm
deg}((z_a^\alpha)^*)=-1$, $a=1,2,\ldots,n$, $\alpha=1,2,\ldots,m$. In this
context one has:
$${\Bbb C}[{\rm Mat}_{mn}]_q^{\rm op}=\bigoplus_{i=0}^\infty{\Bbb C}[{\rm
Mat}_{mn}]_{q,i}^{\rm op},\qquad{\Bbb C}[\overline{\rm
Mat}_{mn}]_q=\bigoplus_{j=0}^\infty{\Bbb C}[\overline{\rm
Mat}_{mn}]_{q,-j}$$
\begin{equation}\label{pka}
{\Bbb C}[{\rm Mat}_{mn}]_q^{\rm op}\otimes{\Bbb C}[\overline{\rm
Mat}_{mn}]_q=\bigoplus_{i,j=0}^\infty{\Bbb C}[{\rm Mat}_{mn}]_{q,i}^{\rm
op}\otimes{\Bbb C}[\overline{\rm Mat}_{mn}]_{q,-j}
\end{equation}
The kernel algebra ${\Bbb C}[[{\rm Mat}_{mn}\times \overline{\rm
Mat}_{mn}]]_q$ will stand for a completion of ${\Bbb C}[{\rm
Mat}_{mn}]_q^{\rm op}\otimes{\Bbb C}[\overline{\rm Mat}_{mn}]_q$ in the
topology associated to the grading in (\ref{pka}). The kernel algebra is
constituted by formal series $\psi=\sum \limits_{i,j=0}^\infty \psi_{ij}$,
with $\psi_{ij}\in{\Bbb C}[{\rm Mat}_{mn}]_{q,i}^{\rm op}\otimes{\Bbb
C}[\overline{\rm Mat}_{mn}]_{q,-j}$.

 Our main result is the following formula for the weighted Bergman kernel:
\begin{equation}\label{wBk}
K_ \lambda=\prod_{j=0}^\infty \left(1+\sum_{i=1}^m(-q^{2(\lambda+j)})^i
\Bbbk_i \right)\cdot \prod_{j=0}^\infty \left(1+\sum_{i=1}^m(-q^{2j})^i
\Bbbk_i \right)^{-1}
\end{equation}
with $\Bbbk_i$ being the polynomial kernels (\ref{pk}). (The right hand side
of (\ref{wBk}) determines an element of ${\Bbb C}[[{\rm Mat}_{mn}\times
\overline{\rm Mat}_{mn}]]_q$ since $\Bbbk_i \in{\Bbb C}[{\rm
Mat}_{mn}]_{q,i}^{\rm op}\otimes{\Bbb C}[\overline{\rm Mat}_{mn}]_{q,-i}$
for all $i=1,2,\ldots,m$).

 In the special case $m=n=1$ we get a well known result \cite{KL}:
$$K_ \lambda=\prod_{j=0}^\infty \left(1-q^{2(\lambda+j)}z \otimes \zeta^*
\right)\cdot \left(\prod_{j=0}^\infty(1-q^{2j}z \otimes \zeta^*)
\right)^{-1}=$$
$$=\sum_{i=0}^\infty \frac{(1-q^{2
\lambda})(1-q^{2(\lambda+1)})\ldots(1-q^{2(\lambda+i-1)})}
{(1-q^2)(1-q^4)\ldots(1-q^{2i})}z^i \otimes \zeta^{*i}.$$
Now passage to a limit as $q \to 1$ and replacement of $\otimes$ by a dot
yields $K_ \lambda \begin{array}{c}\\ \to \\{q \to 1 \atop} \end{array}(1-z
\zeta^*)^{-\lambda}$.

 A q-analogue of an ordinary Bergman kernel for the matrix ball (see
\cite{Hua}) is derivable from (\ref{wBk}) by a substitution $\lambda=m+n$:
$$K=\prod_{j=0}^{m+n-1}\left(1+\sum_{i=1}^m(-q^{2j})^i \Bbbk_i
\right)^{-1}\begin{array}{c}\\ \to \\{q \to 1 \atop} \end{array}\left({\rm
det}(1-{\bf z}\cdot \hbox{\boldmath$\zeta$}^*)\right)^{-(m+n)}.$$

\end{document}